\documentclass[12pt]{article}
\usepackage{amsfonts}
\usepackage{graphicx}
\usepackage{amsmath}
\usepackage{amsthm}
\usepackage{latexsym,bm}
\usepackage{amssymb}
\usepackage{indentfirst}
\begin{document}
\title{{Division Theorems for the Koszul Complex}
\footnotetext{{ Mathematics Classification Primary(2000)}:{
32A38;32W05}\\
\hspace*{7mm}{Keywords: Koszul complex; Skoda triple;  Twisted
version of Skoda's estimate.}}}
\author{Qingchun Ji}

\date{}
\maketitle
\begin{abstract}

We establish a twisted version of Skoda's estimate for the Koszul
complex from which we get global division theorems for the Koszul
complex. This generalizes Skoda's division theorem. We also show how
to use Skoda triples to produce division theorems for the Koszul
complex.

\end{abstract}

\section{Introduction}
Skoda's division theorem is a $L^2$ version of the Corona theorem in
higher dimensions. It turns out to be an important tool in
establishing effective results in complex geometry and algebraic
geometry(see refs.[B87],[El99],[Siu98] and [Siu00]).

Many generalizations of Skoda's division have bee known since [S72].
In [S78] and [D82], division theorems were proved for generically
surjective homomorphisms between holomorphic vector bundles.
Inspired by the Ohsawa-Takegoshi technique, Varolin ([V08]) proved
the twisted version of Skoda's estimate and also introduced Skoda
triple which enabled he to get a series of Skoda-type theorems. By
using the method of residue currents, Andersson also studied
division problem for the Koszul complex and its geometric
applications(see [A04],[AG11] and references therein).

In this paper, we first prove a Skoda-type estimate(see lemma 2) for
the Koszul complex and then try to introduce twisting into such an
apriori estimate. To prove lemma 1, we make use of a generalization
of Skoda's inequality whose proof is included in the appendix. Based
upon the twisted version of our Skoda-type estimate(see lemma 3), we
obtain Skoda-type division theorems for the Koszul complex. In
principal, by using the Skoda triple introduced by Varolin, we can
obtain many examples of division theorems for the Koszul complex.
Moreover, the technique of denominators ([MV07]) could be used to
produce Skoda triples as shown in [V08]. Our main results on the
division problem are theorems 1 and 2. As an application of theorem
1, we give a sufficient condition(see corollary 1) under which the
Koszul complex induces an exact sequence at the level of global
sections. We give explicit examples(corollaries 2 and 3) from
theorem 2 and choosing Skoda triples, then we also discuss the
relations among these results.

\section{Twisted Estimate for the Koszul Complex}

Let $\Omega$  be a domain in\emph{ $\mathbb{C}^{n},$
}$g_{1}\cdots,g_{p}\in\mathcal{O}(\Omega)$, we denote by $g$ the
vector-valued function $(g_{1}\cdots,g_{p})$ where
$\mathcal{O}(\Omega)$ is the ring of holomorphic functions on
$\Omega$. We also denote by $\mathcal{O}_\Omega$ the sheaf of germs
of holomorphic functions on $\Omega$. Now we can introduce the
well-known Koszul complex$$0\rightarrow
{\bigwedge}^p\mathcal{O}_\Omega^{\oplus
p}\overset{\iota_g}{\rightarrow}\cdots\overset{\iota_g}{\rightarrow}{\bigwedge}^\ell\mathcal{O}_\Omega^{\oplus
p}\overset{\iota_g}{\rightarrow}{\bigwedge}^{\ell-1}\mathcal{O}_\Omega^{\oplus
p}\overset{\iota_g}{\rightarrow}\cdots\overset{\iota_g}{\rightarrow}\mathcal{O}
_\Omega\rightarrow 0\eqno{(1)}$$ The sheaf-homomorphism $\iota_g$ is
defined for each $1\leq \ell\leq p$ as follows$$(\iota_g
v)_{i_1\cdots i_{\ell-1}}=\sum_{1\leq i\leq p}g_iv_{ii_1\cdots
i_{\ell-1}},\ \ 1\leq i_1,\cdots, i_{\ell-1}\leq p\eqno{(2)}$$where
$v=(v_{i_1\cdots i_\ell})^p_{i_1\cdots
i_\ell=1}\in{\bigwedge}^\ell\mathcal{O}_\Omega^{\oplus p}$
i.e.$v_{i_1\cdots i_\ell}\in\mathcal{O}_\Omega$ and $v_{i_1\cdots
i_\ell}$ are skew-symmetric on the indices $i_1,\cdots, i_\ell.$\\

We consider, in this paper, the global division problem for the
Koszul complex, i.e. to find sufficient condition of integrability
under which an element $f\in
{\bigwedge}^{\ell-1}\mathcal{O}(\Omega)^{\oplus p}(1\leq\ell\leq p)$
with $\iota_gf=0$ should be contained in the image of $\iota_g,$
that is, there exists some $u\in
{\bigwedge}^\ell\mathcal{O}(\Omega)^{\oplus p}$ such that
$f=\iota_gu$. The celebrated Skoda division theorem is concerned
with the
case where $\ell=1$ of this problem.\\

We will agree on the following index ranges:$$1\leq i,j,k,\ell \leq
p, \ \ \ \ \ 1\leq \alpha,\beta\leq n.$$

In the remaining part of this section, we assume that $\Omega$ is a
bounded domain with smooth boundary and
$g_i\in\mathcal{O}(\Omega)\cap C^\infty(\bar\Omega),1\leq i\leq p,$
which have no common zeros on $\bar\Omega$. Such assumptions will be
drooped in section 2 by standard argument. For functions
$\phi_1,\phi_2\in C^2(\bar\Omega)$, we
define$$\iota_g^*h=e^{\phi_1-\phi_2}\bar{g}\wedge h\eqno{(3)}$$where
$h\in \bigwedge^{\ell-1}\mathcal{O}(\Omega)^{\oplus p},1\leq\ell\leq
p$ and $\bar{g}\wedge h$ is given by $$(\bar{g}\wedge
h)_{i_1\cdots\i_{\ell}}=-\sum_{1\leq\sigma\leq\ell}(-1)^{\sigma}\bar{g}_{i_\sigma}h_{i_1\cdots\hat{i}_\sigma\cdots
i_{\ell}}.\eqno{(4)}$$ where $\hat{i}_\sigma$ means that the index
$i_\sigma$ is omitted. To formulate our a priori estimate, we need
to introduce the following space.$$F=\{h\in
{\bigwedge}^{\ell-1}\mathcal{O}(\Omega)^{\oplus p} \mid \iota_g
h=0\}.\eqno{(5)}$$

Now we start to estimate
$\|\iota_g^*h+\bar{\partial}_{\phi_1}^*v\|_{\phi_1}^2$ from below
for an arbitrary $h\in F$ which is smooth on $\bar\Omega$ and $v\in
{\rm Dom}\bar{\partial}_{\phi_1}^*\subseteq{\bigwedge}^\ell
L^2_{0,1}(\Omega,\phi_1)^{\oplus p}$ satisfying $\bar\partial v=0$.
The reason why we should estimate such a term lies in Skoda's
fundamental lemma(see lemma 4 at the end of this section).

\begin{eqnarray*}\|\iota_g^*h+\bar{\partial}_{\phi_1}^*v\|_{\phi_1}^2 &=&
\|\iota_g^*h\|_{\phi_1}^2+2{\rm Re}(\iota_g^*h,\bar{\partial}_{\phi_1}^*v)_{\phi_1}+\|\bar{\partial}_{\phi_1}^*v\|_{\phi_1}^2\\
&=& I+II+III.\end{eqnarray*}

By definition, we have
$$I=\int_\Omega\|\bar{g}\wedge
h\|^2e^{\phi_1-2\phi_2}dV.$$ Form the following identity
\begin{eqnarray*}\|\bar{g}\wedge h\|^2&=&\frac{1}{\ell
!}\sum_{\overset{1\leq\sigma,\delta\leq \ell}{1\leq i_1,\cdots,
i_\ell\leq
p}}(-1)^{\sigma+\delta}\bar{g}_{i_\sigma}g_{i_\delta}h_{i_1\cdots\hat{i}_\sigma\cdots
i_\ell}\bar{h}_{i_1\cdots\hat{i}_\delta\cdots i_\ell}\\
&=&\frac{1}{\ell !}\sum_{\overset{1\leq\sigma\leq \ell}{1\leq
i_1,\cdots, i_\ell\leq
p}}|g_{i_\sigma}|^2|h_{i_1\cdots\hat{i}_\sigma\cdots i_\ell}|^2\\
&\ &-\frac{1}{\ell !}\sum_{\overset{1\leq\sigma\neq\delta\leq
\ell}{1\leq i_1,\cdots,\hat{i}_\sigma,\cdots,\hat{i}_\delta,\cdots,
i_\ell\leq p}}\sum_{1\leq i_\sigma\leq
p}\bar{g}_{i_\sigma}\bar{h}_{i_\sigma
i_1\cdots\hat{i}_\sigma\cdots\hat{i}_\delta\cdots i_\ell}\sum_{1\leq
i_\delta\leq p}g_{i_\delta}h_{i_\delta
i_1\cdots\hat{i}_\sigma\cdots\hat{i}_\delta\cdots i_\ell}\\
&\overset{\iota_gh=0}{=}& \frac{1}{(\ell -1)!}\sum_{1\leq i,
i_1,\cdots, i_{\ell-1}\leq
p}|g_{i}|^2|h_{i_1\cdots i_{\ell-1}}|^2\\
&=&\|g\|^2\|h\|^2,\end{eqnarray*}it follows
that$$I=\int_\Omega\|g\|^2\|h\|^2e^{\phi_1-2\phi_2}dV$$where
$\|g\|^2=\displaystyle{\sum_{1\leq i\leq
p}}|g_i|^2,\|h\|^2=\frac{1}{(\ell-1)!}\displaystyle{\sum_{1\leq
i_1,\cdots, i_{\ell-1}\leq p}}|h_{i_1\cdots i_{\ell-1}}|^2.$

Assuming that $\phi_1,\phi_2$ are related as follows
$$\phi_2=\phi_1+log\|g\|^2,\eqno{(6)}
$$then we have$$I=\|h \|_{\phi_2}^2.\eqno{(7)}$$
Moreover, with $\phi_1,\phi_2$ chosen in such a way, we
also
have$$\iota_g^*=\frac{\bar{g}}{\|g\|^2}\wedge.$$\\

To deal with the term $II$, we first use integration by parts to
move the operator $\bar{\partial}_{\phi_1}$ from right to
left.\begin{eqnarray*}2{\rm
Re}(\iota_g^*h,\bar{\partial}_{\phi_1}^*v)_{\phi_1} &=& 2{\rm
Re}\int_\Omega(\frac{\bar{g}}{\|g\|^2}\wedge
h,\bar{\partial}_{\phi_1}^*v)e^{-\phi_1}dV \\
&=& 2{\rm
Re}\int_\Omega(\bar{\partial}(\frac{\bar{g}}{\|g\|^2}\wedge
h),v)e^{-\phi_1}dV\\ &=& 2{\rm Re}\int_\Omega(\sum_{1\leq\alpha\leq
n}\overline{\partial_\alpha(\frac{g}{\|g\|^2})}\wedge hd\bar
z_\alpha,v)e^{-\phi_1}dV.
\end{eqnarray*}Let
$$v_{i_1\cdots i_{\ell}}=\sum_{1\leq \alpha\leq n}v_{i_1\cdots
i_{\ell}\alpha}d\bar{z}_\alpha,$$then we have
\begin{eqnarray*}II&=& \frac{-2{\rm
Re}}{\ell!}\int_\Omega
\sum_{\overset{1\leq\sigma\leq\ell}{\overset{1\leq\alpha\leq
n}{1\leq i_1,\cdots,i_{\ell}\leq p}}}(-1)^\sigma
\overline{\partial_\alpha(\|g\|^{-2}g_{i_\sigma})}h_{i_1\cdots\hat{i}_\sigma\cdots
i_{\ell}}\overline{v_{i_1\cdots i_{\ell}\alpha}}e^{-\phi_1}dV\\&=&
\frac{2{\rm Re}}{(\ell-1)!}\int_\Omega\sum_{\overset{1\leq\alpha\leq
n}{1\leq i,i_1,\cdots,i_{\ell-1}\leq p}}h_{i_1\cdots
i_{\ell-1}}\overline{\partial_\alpha(\|g\|^{-2}g_i)v_{ii_1\cdots
i_{\ell-1}\alpha}}e^{-\phi_1}dV\end{eqnarray*} which
gives\begin{eqnarray*}II &\geq&
-\int_\Omega\frac{1}{b}\|h\|^2e^{-\phi_2}dV \\& \
&-\frac{1}{(\ell-1)!}\sum_{1\leq i_1,\cdots,i_{\ell-1}\leq
p}\int_\Omega b\|g\|^2\left|\sum_{\overset{1\leq\alpha\leq n}{1\leq
i\leq p}}\partial_\alpha(\|g\|^{-2}g_i)v_{ii_1\cdots
i_{\ell-1}\alpha}\right|^2e^{-\phi_1}dV\end{eqnarray*} where $b>1$
is an arbitrary measurable function on $\Omega$.\\

To handle the second term in the above inequality, we need the
following lemma which is a generalization of Skoda's inequality([S72], page 552).\\

{\bf Lemma 1} Given constants $a_i,b_{i\alpha},c_{i_1\cdots
i_\ell\alpha}\in \mathbb{C}(1\leq i,\ell,i_1,\cdots,i_\ell\leq
p,1\leq \alpha\leq n,p,n\in \mathbb{N})$  where $c_{i_1\cdots
i_\ell\alpha}$
 are skew-symmetric with respect to $i_1\cdots
i_\ell$. Then for any $1\leq i_1<\cdots<i_{\ell-1}\leq p$, it holds
that$$\left|\sum_{\overset{1\leq i,j\leq p}{1\leq \alpha\leq
n}}\overline{a_j}(a_jb_{i\alpha}-a_ib_{j\alpha})c_{ii_1\cdots
i_{\ell-1}\alpha}\right|^2\leq q\sum_{1\leq i\leq
p}|a_i|^2\sum_{\overset{1\leq i\leq p}{1\leq j<k\leq
p}}\left|\sum_{1\leq\alpha\leq n}(a_jb_{k\alpha}-a_kb_{j\alpha})c_{i
i_1\cdots i_{\ell-1}\alpha}\right|^2\eqno{(8)}$$ where $q$ is the
positive constant defined by
\begin{equation*}q=\begin{cases}
\min\{p-1,n\}, & \ell=1;\\\min\{p-\ell+1,n\}, &\ell\geq 2.\end{cases}\end{equation*}\\

We postpone our proof of lemma 1 to the appendix. Now we continue
the estimate for $II$. Applying lemma 1 with $$a_i=g_i,\ \
b_{i\alpha}=\partial_\alpha g_i,\ \ c_{i_1\cdots
i_\ell}=v_{i_1\cdots i_\ell},$$ we obtain for fixed $1\leq
i_1<\cdots<i_{\ell-1}\leq p$
that\begin{eqnarray*}\|g\|^2\left|\sum_{\overset{1\leq\alpha\leq
n}{1\leq i\leq p}}\partial_\alpha(\|g\|^{-2}g_i)v_{ii_1\cdots
i_{\ell-1}\alpha}\right|^2 &=&
\|g\|^{-6}\left|\sum_{\overset{1\leq\alpha\leq n}{1\leq i,j\leq
p}}\overline{g_j}(g_j\partial_\alpha g_i-g_i\partial_\alpha
g_j)v_{ii_1\cdots i_{\ell-1}\alpha}\right|^2 \\
&\leq& q\|g\|^{-4}\sum_{\overset{1\leq i\leq p}{1\leq j<k\leq
p}}\left|\sum_{1\leq\alpha\leq n}(g_j\partial_\alpha
g_k-g_k\partial_\alpha g_j)v_{ii_1\cdots i_{\ell-1}\alpha}\right|^2
\\ &\overset{(*)}{=}& q\|g\|^{-2}\sum_{1\leq j,k\leq p}\left|\sum_{1\leq\alpha\leq n}
\partial_\alpha g_kv_{j i_1\cdots i_{\ell-1}\alpha}\right|^2\\ &\ &-q\|g\|^{-4}\sum_{1\leq j\leq p}
\left|\sum_{\overset{1\leq k\leq p}{1\leq \alpha\leq
n}}g_k\partial_\alpha g_kv_{j i_1\cdots i_{\ell-1}\alpha}\right|^2 \\
&=& q\sum_{\overset{1\leq j\leq p}{1\leq\alpha,\beta\leq n}}\
\partial_\alpha\partial_{\bar{\beta}}\log\|g\|^2\cdot v_{j i_1\cdots
i_{\ell-1}\alpha}\overline{v_{j i_1\cdots
i_{\ell-1}\beta}}\end{eqnarray*}where $q$ is the constant in lemma
1. We have used the Lagrange identity to get $(*)$.

Consequently, we can estimate $II$ from below as follows.$$
II\geq-\int_\Omega\frac{1}{b}\|h\|^2e^{-\phi_2}dV \ \ \ \ \ \ \ \ \
\ \ \ \ \ \ \ \ \ \ \ \ \ \ \ \ \ \ \ \ \ \ \ \ \ \ \ \ \ \ \ \ \ \
\ \ \ \ \ \ \ $$$$\ \ \ \ \ \ \ \ \
-\frac{q}{(\ell-1)!}\sum_{\overset{1\leq i_1,\cdots,i_{\ell}\leq
p}{1\leq\alpha,\beta\leq n}}\int_\Omega
b\partial_\alpha\partial_{\bar{\beta}}\log\|g\|^2\cdot v_{i_1\cdots
i_{\ell}\alpha}\overline{v_{i_1\cdots i_{\ell}\beta}}e^{-\phi_1}dV. \eqno{(9)}$$\\

Since $I$ only involves $h$, we want to control the second term in
(9) by $III$. By using the standard
Kohn-Morrey-H$\rm{\ddot{o}}$rmander identity and the boundary
condition $v\in {\rm Dom}\bar\partial_{\phi_1}^*,$ we can estimate,
in the case where $\Omega$ is assumed additionally to be
pseudoconvex, the last term $III$ as follows

$$III=\frac{1}{\ell!}\sum_{1\leq i_1,\cdots,i_{\ell}\leq
p}\int_\Omega\left|\bar{\partial}^*_{\phi_1}v_{i_1\cdots
i_\ell}\right|^2e^{-\phi_1}dV\ \ \ \ \ \ \ \ \ \ \ \ \ \ \ \ \ \ \ \
\ \ \ \ \ \ \ \ \ \
$$
$$
\geq\frac{1}{\ell!}\sum_{\overset{1\leq i_1,\cdots,i_{\ell}\leq
p}{1\leq\alpha,\beta\leq
n}}\int_\Omega\partial_\alpha\partial_{\bar{\beta}}\phi_1v_{i_1\cdots
i_{\ell}\alpha}\overline{v_{i_1\cdots i_{\ell}\beta}}e^{-\phi_1}dV.\
\ \ \ \ \ \ \ \ \ \ \ \ \eqno{(10)}$$ for any $v\in {\rm
Dom}\bar{\partial}_{\phi_1}^*\subseteq{\bigwedge}^\ell
L^2_{0,1}(\Omega,\phi_1)^{\oplus p}$ satisfying $\bar\partial v=0$.\\

Taking the sum of (7),(9) and (10), we get the desired Skoda-type
estimate for the Koszul complex.\\

{\bf Lemma 2} Let $\Omega$ be a bounded pseudoconvex domain with
smooth boundary and $g_i\in\mathcal{O}(\Omega)\cap
C^\infty(\bar\Omega)(1\leq i\leq p)$ which have no common zeros on
$\bar{\Omega}$. We assume that $\phi_1,\phi_2\in C^2(\bar{\Omega})$
are functions satisfying (6), and $b>1$ is a measurable function on
$\Omega$. Then for any $h\in F$ and any $v\in {\rm
Dom}\bar\partial_{\phi_1}^*\subseteq{\bigwedge}^\ell
L^2_{0,1}(\Omega,\phi_1)^{\oplus p}$ satisfying $\bar\partial v=0,$
it holds that

$$\|\iota_g^*h+\bar{\partial}_{\phi_1}^*v\|_{\phi_1}^2 \geq
\int_\Omega\frac{b-1}{b}\|h\|^2e^{-\phi_2}dV \ \ \ \ \ \ \ \ \ \ \ \
\ \ \ \ \ \ \ \ \ \ \ \ \ \ \ \ \ \ \ \ \ \ \ \ \ \ \ \ \ \ \ \ \ \
\ \ \ \ \ \
$$$$ \ \ \ \ \ \ \ \ \ \ \ \ \ \ \ \ \ \ \ \ \ \ \ \ \ \ \ \ \ \ \
+\frac{1}{\ell!}\sum_{\overset{1\leq i_1,\cdots,i_{\ell}\leq
p}{1\leq\alpha,\beta\leq
n}}\int_\Omega(\partial_\alpha\partial_{\bar{\beta}}\phi_1 -q \ell
b\partial_\alpha\partial_{\bar{\beta}}\log\|g\|^2 ) v_{i_1\cdots
i_{\ell}\alpha}\overline{v_{i_1\cdots i_{\ell}\beta}}e^{-\phi_1}dV
.\eqno{(11)}$$\\

Now we want to introduce twisting into the apriori estimate (11).
Following [V08], we twist simultaneously the weights $\phi_1,\phi_2$
by a function $0<a\in C^2(\bar\Omega)$ and consider the following
new weights.$$\varphi_1=\phi_1+\log a,\ \ \ \ \varphi_2=\phi_2+\log
a.\eqno{(12)}$$From (6) and (12), it follows that
$$\varphi_2=\varphi_1+log\|g\|^2.\eqno{(13)}$$ By the definition (3) of $\iota_g^*$,
we know
$$\iota_g^*=e^{\varphi_1-\varphi_2}\bar{g}\wedge.\eqno{(14)}$$
From the definition of $\bar\partial^*_{\phi_1}$ and (12), we get
$$(\bar\partial^*_{\phi_1}v)_{i_1\cdots
i_{\ell}} = (\bar\partial^*_{\varphi_1}v)_{i_1\cdots
i_{\ell}}-\sum_{1\leq\alpha\leq n}\frac{\partial_\alpha
a}{a}v_{i_1\cdots i_{\ell}\alpha}.$$which implies the following
identity
\begin{eqnarray*}\|\sqrt a\iota_g^*h+\sqrt
a\bar\partial_{\varphi_1}^*v\|^2_{\varphi_1}
&=&\|\iota_g^*h+\bar\partial_{\varphi_1}^*v\|^2_{\phi_1}\\ &=&
\|\iota_g^*h+\bar\partial_{\phi_1}^*v\|^2_{\phi_1} \\
&\ &+\frac{1}{\ell !} \sum_{\overset{1\leq \alpha,\beta\leq 1}{1\leq
i_1,\cdots, i_\ell\leq p}}\int_\Omega a^{-2}\partial_\alpha
a\partial_{\bar\beta}av_{i_1\cdots i_\ell\alpha}
\overline{v_{i_1\cdots i_\ell\beta}}e^{-\phi_1}dV\\
&\ & + \frac{2}{\ell !}{\rm Re}\sum_{1\leq i_1,\cdots, i_\ell\leq
p}\int_\Omega
\frac{1}{a}(\iota_g^*h+\bar\partial_{\phi_1}^*v)_{i_1\cdots
i_\ell}\overline{\sum_{1\leq\alpha\leq n}\partial_\alpha
av_{i_1\cdots i_\ell\alpha}}e^{-\phi_1}dV
\\ &=&
\|\iota_g^*h+\bar\partial_{\phi_1}^*v\|^2_{\phi_1} \\
&\ &-\frac{1}{\ell !} \sum_{\overset{1\leq \alpha,\beta\leq 1}{1\leq
i_1,\cdots, i_\ell\leq p}}\int_\Omega a^{-2}\partial_\alpha
a\partial_{\bar\beta}av_{i_1\cdots i_\ell\alpha}
\overline{v_{i_1\cdots i_\ell\beta}}e^{-\phi_1}dV\\
&\ & + \frac{2}{\ell !}{\rm Re}\sum_{1\leq i_1,\cdots, i_\ell\leq
p}\int_\Omega\frac{1}{a}(\iota_g^*h+\bar\partial_{\varphi_1}^*v)_{i_1\cdots
i_\ell}\overline{\sum_{1\leq\alpha\leq n}\partial_\alpha
av_{i_1\cdots i_\ell\alpha}}e^{-\phi_1}dV.\end{eqnarray*}
Substituting (11) and the following equation
$$\partial_\alpha\partial_{\bar\beta}\phi_1 =
\partial_\alpha\partial_{\bar\beta}\varphi_1 -
a^{-1}\partial_\alpha\partial_{\bar\beta}a + a^{-2}\partial_\alpha
a\partial_{\bar\beta}a$$ into the above identity, it follows that
\begin{eqnarray*}\|\sqrt a\iota_g^*h+\sqrt
a\bar\partial_{\varphi_1}^*v\|^2_{\varphi_1} &\geq&
\int_\Omega\frac{b-1}{b}\|h\|^2e^{-\phi_2}dV \\
&\ &+\frac{1}{\ell !} \sum_{\overset{1\leq \alpha,\beta\leq 1}{1\leq
i_1,\cdots, i_\ell\leq p}}\int_\Omega [\partial_\alpha
\partial_{\bar\beta}\varphi_1-q\ell b\partial_\alpha
\partial_{\bar\beta}\log\|g\|^2\\
 &\ &\ \ \ \ \ \ \ \ \ \ \ \ \ \ \ \ \  -\frac{1}{a}\partial_\alpha
\partial_{\bar\beta}a-\frac{1}{a\lambda}\partial_\alpha a
\partial_{\bar\beta}a] v_{i_1\cdots i_\ell\alpha}
\overline{v_{i_1\cdots i_\ell\beta}}e^{-\phi_1}dV\\
&\ & -\|\sqrt
\lambda\iota_g^*h+\sqrt\lambda\bar\partial_{\varphi_1}^*v\|^2_{\varphi_1}.
\end{eqnarray*}where $\lambda>0$ is a measurable function on
$\Omega$.

If the following condition holds
$$a\partial_\alpha\partial_{\bar{\beta}}\phi_1-\partial_\alpha\partial_{\bar{\beta}}a
-\lambda^{-1}\partial_\alpha a\partial_{\bar{\beta}}a \geq q\ell
ab\partial_\alpha\partial_{\bar{\beta}}\log\|g\|^2\eqno{(15)}$$
where both sides are understood as symmetric sesquilinear forms and
$q$ is the constant in lemma 1, then we have$$\|\sqrt{a+
\lambda}\iota_g^*h+\sqrt{a+\lambda}\bar\partial_{\varphi_1}^*v\|^2_{\varphi_1}
\geq
\int_\Omega\frac{(b-1)a}{b}\|h\|^2e^{-\varphi_2}dV. \eqno{(16)}$$  where $\iota_g^*$ is given by (14). The estimate (16) is a twisted version of (11).\\

We summarize previously obtained estimates in the following lemma.\\

{\bf Lemma 3} Let $\Omega$ be a bounded pseudoconvex domain with
smooth boundary and $g_i\in\mathcal{O}(\Omega)\cap
C^\infty(\bar\Omega)(1\leq i\leq p)$ which have no common zeros on
$\bar{(\Omega)}$. We assume that $\varphi_1,\varphi_2\in
C^2(\bar{\Omega})$ are functions satisfying (13), $0<a\in
C^2(\bar{\Omega})$ and $1<b,0<\lambda$ are measurable functions on
$\Omega$. Then for any $h\in F$ and any $v\in {\rm
Dom}\bar\partial_{\varphi_1}^*\subseteq{\bigwedge}^\ell
L^2_{0,1}(\Omega,\varphi_1)^{\oplus
p}$ satisfying $\bar\partial v=0,$ the twisted estimate (16) follows from the condition (15).\\

The next fundamental lemma reduces the  problem of establishing
division theorems to an apriori estimate(see [S72] and [V08] for proofs).\\

{\bf Lemma 4} Let $H,H_{0},H_{1},H_{2}$ be Hilbert spaces,
$T:H_{0}\rightarrow H$ be a bounded linear operator,
$T_{\nu}:H_{\nu-1}\rightarrow H_{\nu}(\nu=1,2)$ be linear, closed,
densely defined operators such that $T_{2}\circ T_{1}=0,$ and let
$F\subseteq H$ be a closed subspace such that $T({\rm
Ker}T_1)\subseteq F.$ Then for every $f\in F$ and constant $C>0$ the
following statements are equivalent

1. There exists at least one $u\in \rm{Ker}T_{1}$ such that $Tu=f$,
$\left\Vert u\right\Vert _{H_0}\leq C;$

2. $|(f,h)_{H}|\leq \left\Vert T^{*}h+T_{1}^{*}v\right\Vert
_{H_{0}}$ holds for any $h\in
F,v\in {\rm Dom}T_{1}^{*}\cap {\rm Ker}T_{2}.$\\

To apply lemma 4, we consider, for any fixed $1\leq \ell\leq p$, the
following Hilbert spaces and
operators.$$H_0={\bigwedge}^{\ell}L^2(\Omega,\varphi_1)^{\oplus p},\
H_1={\bigwedge}^{\ell}L_{0,1}^2(\Omega,\varphi_1)^{\oplus
p},$$$$H_2={\bigwedge}^{\ell}L_{0,2}^2(\Omega,\varphi_1)^{\oplus
p},H={\bigwedge}^{\ell-1}L^2(\Omega,\varphi_2)^{\oplus
p},$$$$T=\sqrt{a+\lambda}\circ\iota_g,
T_1=\bar\partial\circ\sqrt{a+\lambda},
T_2=\sqrt{a}\circ\bar\partial.$$where these functions
$\varphi_1,\varphi_2,0<a,0<\lambda\in
C^2(\bar\Omega)$ will be determined later.\\

Since $\Omega$ is assumed to be bounded and $a,\lambda \in
C^2(\bar\Omega)$, the operator $T$ is a bounded linear mapping from
$H_0$ to $H$. $T_1,T_2$ are, by definition, both densely defined and
closed. The space $F$ defined by (5) is obviously a closed subspace
of $H$. From the definition of $T_1$, we have
$${\rm Ker}T_1\subseteq\frac{1}{\sqrt{a+\lambda}}{\bigwedge}^\ell\mathcal{O}(\Omega)^{\oplus
p}.$$Consequently,  $\iota_g^2=0$ implies that
$$F\supseteq T({\rm Ker}T_1).$$\\ It is also easy to see that the adjoint of $T$ and $T_1$ are given by
$$T^*=\sqrt{a+\lambda}\iota_g^*,\ \ T_1^*=\sqrt{a+\lambda}\bar{\partial}_{\varphi_1}^*$$where $\iota^*_g$ is the mapping in (14).\\

\section{Division Theorems}

First we apply lemma 3 in the simplest situation where the function
$a$ is a constant to establish a division theorem for the Koszul
complex. We denote by ${\rm PSH}(\Omega)$
the set of plurisubharmonic functions on $\Omega$.\\

{\bf Theorem 1} Let $\Omega\subseteq\mathbb{C}^n$ be a pseudoconvex
domain, $g_i\in\mathcal{O}(\Omega)(1\leq i\leq p)$, $\psi\in{\rm
PSH}(\Omega)$ and $\tau >1$ be a constant. For every $f\in
{\bigwedge}^{\ell-1}\mathcal{O}(\Omega)^{\oplus p}$, if $\iota_gf=0$
and
$$\int_\Omega\|f\|^2\|g\|^{-2(q\ell\tau+1)}e^{-\psi}dV<\infty,\eqno{(17)}$$then
there exists an $u\in {\bigwedge}^{\ell}\mathcal{O}(\Omega)^{\oplus
p}$ such that
$$\iota_gu=f,\ \ \int_\Omega\|u\|^2\|g\|^{-2q\ell\tau}e^{-\psi}dV
\leq\frac{\tau}{\tau-1}\int_\Omega\|f\|^2\|g\|^{-2(q\ell\tau+1)}e^{-\psi}dV\eqno{(18)}$$where
$p\in\mathbb{N}, 1\leq\ell\leq p$ and $q$ is the constant in lemma 1.\\

Proof. By the standard argument of smooth approximation, the
holomorphic extension technique and taking weak limit(proceed as
[S72] and [D82]), we can assume without loss of generality that
$\Omega$ is a bounded pseudoconvex domain with smooth boundary,
$g_i\in\mathcal{O}(\Omega)\cap C^\infty(\bar\Omega)(1\leq i\leq p)$
have no common zeros on $\bar{\Omega}$ and $\psi\in{\rm
PSH}(\Omega)\cap C^\infty(\bar\Omega)$.

Given a constant $\tau>1$, we can always find constants
$0<\lambda<1<b$ such that
$$\tau=\frac{b}{1-\lambda}.$$Set$$a=1-\lambda,$$ then the functions
$$\varphi_1=q\ell\tau \log\|g\|^2 +\psi,\ \ \varphi_2=(q\ell\tau+1)\log\|g\|^2
+\psi$$satisfy the conditions (13) and (15). In this case, we have
$T=\iota_g$ and $T_1=\bar\partial.$

Let $F$ be the closed subspace defined by (5) in section 2 and $h\in
F$, then we get by lemma 3 that\begin{eqnarray*}|(f,h)_{H}|^2 &\leq&
\int_\Omega\frac{(b-1)(1-\lambda)}{b}\|h\|^2e^{-\varphi_2}dV\int_\Omega\frac{b}{(b-1)(1-\lambda)}\|f\|^2e^{-\varphi_2}dV \\
&\leq&\int_\Omega\frac{b}{(b-1)(1-\lambda)}\|f\|^2e^{-\varphi_2}dV\cdot\|T^*h+T_1^*v\|_{H_0}^2.\end{eqnarray*}
It follows from lemma 4 that there is an $u_\lambda\in H_0$ such
that
$$f=Tu_\lambda=\iota_g u_\lambda$$ and the weighted $L^2$ norm of $u$ could be
estimated as follows
\begin{eqnarray*}\int_\Omega\|u_\lambda\|^2\|g\|^{-2q\ell\tau}e^{-\psi}dV&=&\|u_\lambda\|^2_{H_0}\\
&\leq&
\int_\Omega\frac{b}{(b-1)(1-\lambda)}\|f\|^2e^{-\varphi_2}dV\\&=&\int_\Omega\frac{\tau}{(1-\lambda)(\tau+\lambda-1)}\|f\|^2e^{-\varphi_2}dV
\\&=&\int_\Omega\frac{\tau}{(1-\lambda)(\tau+\lambda-1)}\|f\|^2\|g\|^{-2(q\ell\tau+1)}e^{-\psi}dV.\end{eqnarray*}
The desired solution $u$ follows from the above inequality by taking
weak limit of $u_\lambda$ as
$\lambda\rightarrow 0 +$.Q.E.D.\\

If $g_1,\cdots,g_p$ have no common zeros, then for any $f\in
{\bigwedge}^{\ell-1}\mathcal{O}(\Omega)^{\oplus p}$, we can use the
plurisubharmonic exhaustion function of $\Omega$ to construct a
plurisubharmonic weight function $\psi$ on $\Omega$ such that
$$\int_{\Omega}\|f\|^2\|g\|^{-2(q\ell\tau+1)}e^{-\psi}dV<\infty.$$ Applying theorem 1 on
$\Omega$, we know there exists some $u\in
{\bigwedge}^\ell\mathcal{O}(\Omega)^{\oplus p}$ such that
$f=\iota_gu$ holds on $\Omega$. This gives the following corollary.\\

{\bf Corollary 1}  Let $\Omega\subseteq\mathbb{C}^n$ be a
pseudoconvex domain, $g_i\in\mathcal{O}(\Omega)(1\leq i\leq p)$. If
$g_1,\cdots,g_p$ have no common zeros, then the Koszul complex (1)
induces an exact sequence at the level of global sections, i.e. for
every $f\in {\bigwedge}^{\ell-1}\mathcal{O}(\Omega)^{\oplus
p}(1\leq\ell\leq p)$ satisfying $\iota_gf=0$ there is some $u\in
{\bigwedge}^\ell\mathcal{O}(\Omega)^{\oplus p}$ such that
$f=\iota_gu$.\\

{\bf Remark 1} (i) The special case of theorem 1 when $\ell=1$ is
exactly the celebrated Skoda division theorem([S72]). If we make use
of lemma 2 instead of lemma 3, the proof of theorem 1 will be a
little bit easier.

(ii) When the common zero locus of $g_1,\cdots,g_p$ is empty, it is
easy to see that the Koszul complex (1) provides a resolution of
$\mathcal{O}_\Omega={\bigwedge}^p\mathcal{O}_\Omega^{\oplus p}$.
Thus corollary 1 also
follows from Cartan's theorem B and the De Rham-Weil isomorphism theorem.\\

To establish division theorems with a nonconstant function $a$ in
(12), we use the technique of Skoda triple which was introduced by
Varolin([V08]). We first recall the definition of a Skoda triple.\\

{\bf Definition } A Skoda triple $(\varphi,F,q)$ consists of a
positive integer $q$ and $C^2$ functions
$\varphi:(1,\infty)\rightarrow\mathbb{R},F:(1,\infty)\rightarrow
\mathbb{R}$ such that
$$x+F(x)>0,[x+F(x)]\varphi^{'}(x)+F^{'}(x)+1>0\ {\rm and}\ [x+F(x)]\varphi^{''}(x)+F^{''}(x)< 0$$
hold for every $x>1.$\\

It is easy to see that $(\varepsilon\log x,0,q)$ is a Skoda triple
where $\varepsilon$ is a positive constant and $q$ is a positive
integer. This example was shown in [V08] where the technique of
denominators was also used to construct Skoda triple of the type
$(0,F,q)$.

Based upon the apriori estimate (16) and lemma 4, the notion of
Skoda triple is quite useful to produce examples of division
theorems.\\

{\bf Theorem 2} Let $\Omega\subseteq\mathbb{C}^n$ be a pseudoconvex
domain, $g_i\in\mathcal{O}(\Omega)(1\leq i\leq p)$, $\psi\in{\rm
PSH}(\Omega)$. We assume that $$\|g\|<1\ \text{holds on }\Omega.$$
For every $f\in {\bigwedge}^{\ell-1}\mathcal{O}(\Omega)^{\oplus p}$,
if $\iota_gf=0$ and
$$\int_\Omega\|f\|^2\frac{b}{a(b-1)}\|g\|^{-2(q\ell+1)}e^{\varphi\circ\xi-\psi}dV<\infty,\eqno{(19)}$$then
there exists an $u\in {\bigwedge}^{\ell}\mathcal{O}(\Omega)^{\oplus
p}$ such that $\iota_gu=f$ and
$$\int_\Omega\|u\|^2\frac{1}{(a+\lambda)}{\|g\|^{-2q\ell}}e^{\varphi\circ\xi
-\psi}dV
\leq\int_\Omega\|f\|^2\frac{b}{a(b-1)}\|g\|^{-2(q\ell+1)}e^{\varphi\circ\xi-\psi}dV\eqno{(20)}$$where
$p\in\mathbb{N}, 1\leq\ell\leq
p,\xi=1-\log\|g\|^2,a=\xi+F\circ\xi,b=\frac{a\varphi^{'}\circ\xi+F^{'}\circ\xi+1}{qa\ell}+1,\lambda=\Lambda\circ\xi,
\Lambda(x)=\frac{-(1+F^{'}(x))^2}{F^{''}(x)+(x+F(x))\varphi^{''}(x)},$
$(\varphi,F,q)$ is a Skoda triple and $q$ is the constant in lemma 1.\\

Proof. Given a Skoda triple, we start to construct functions
$\varphi_1,\varphi_2,a>0,\lambda>0$ and $b>1$ which satisfy
conditions (13) and (15).

Set
\begin{eqnarray*}\varphi_1&=&-\varphi\circ\xi+\psi+q\ell\log\|g\|^2,\\
\varphi_2&=&-\varphi\circ\xi+\psi+(q\ell+1)\log\|g\|^2\end{eqnarray*}
then we
get\begin{eqnarray*}a\partial_\alpha\partial_{\bar{\beta}}\varphi_1-\partial_\alpha\partial_{\bar{\beta}}a
-\lambda^{-1}\partial_\alpha a\partial_{\bar{\beta}}a &=&
(a\varphi^{'}\circ\xi+F^{'}\circ\xi+1+q\ell
a)\partial_\alpha\partial_{\bar\beta}log\|g\|^2\\
&\ &
-[a\varphi^{''}\circ\xi+F^{''}\circ\xi+\lambda^{-1}(1+F^{'}\circ\xi)^2]
\partial_\alpha\xi\partial_{\bar\beta}\xi\\
&=&(a\varphi^{'}\circ\xi+F^{'}\circ\xi+1+q \ell
a)\partial_\alpha\partial_{\bar\beta}log\|g\|^2.\end{eqnarray*}The
last equality follows from the definition of $\lambda.$

Now it suffices to choose $b>1$ such that
$$a\varphi^{'}\circ\xi+F^{'}\circ\xi+1=qa\ell(b-1),$$
i.e.$$b=\frac{a\varphi^{'}\circ\xi+F^{'}\circ\xi+1}{qa\ell}+1.$$ By
repeating the argument in the proof of theorem 1, we obtain some
$\tilde {u}\in H_0={\bigwedge}^{\ell}L^2(\Omega,\varphi_1)^{\oplus
p}$ satisfying $$\bar\partial\sqrt{a+\lambda}\tilde{u}=0,\ \ \
\int_\Omega\|\tilde{u}\|^2{\|g\|^{-2q\ell}}e^{\varphi\circ\xi
-\psi}dV
\leq\int_\Omega\|f\|^2\frac{b}{a(b-1)}\|g\|^{-2(q\ell+1)}e^{\varphi\circ\xi-\psi}dV.$$
Thus we get the desired solution $u=\sqrt{a+\lambda}\tilde{u}$.Q.E.D.\\

If we take into account the special Skoda triple $(\varepsilon\log
x, 0, q)$ where $\varepsilon$ is a positive constant and $q$ is the
constant in lemma 1, applying theorem 2 to this Skoda triple, we have the following corollary.\\

{\bf Corollary 2} Let $\Omega\subseteq\mathbb{C}^n$ be a
pseudoconvex domain, $g_i\in\mathcal{O}(\Omega)(1\leq i\leq p)$,
$\psi\in{\rm PSH}(\Omega)$. We assume that $$\|g\|<1\ \text{holds on
}\Omega.$$ For every $f\in
{\bigwedge}^{\ell-1}\mathcal{O}(\Omega)^{\oplus p}$, if $\iota_gf=0$
and
$$\int_\Omega\|f\|^2
\frac{(1-\log\|g\|^2)^\varepsilon}{\|g\|^{2(q\ell+1)}}
e^{-\psi}dV<\infty,\eqno{(21)}$$then there exists some $u\in
{\bigwedge}^{\ell}\mathcal{O}(\Omega)^{\oplus p}$ such that
$\iota_gu=f$ and
$$\int_\Omega\|u\|^2\frac{(1-\log\|g\|^2)^{\varepsilon -1}}{ \|g\|^{2q\ell}}e^{-\psi}dV
\leq\frac{q\ell+\varepsilon+1}{\varepsilon}\int_\Omega\|f\|^2\frac{(1-\log\|g\|^2)^\varepsilon}{\|g\|^{2(q\ell+1)}}
e^{-\psi}dV\eqno{(22)}$$where $p\in\mathbb{N}, 1\leq\ell\leq p,$
$\varepsilon>o$ is a constant and $q$ is the constant in lemma 1.\\

Proof. For the given Skoda triple $(\varepsilon\log x, 0, q)$, we
have
$$[x+F(x)]\varphi^{'}(x)+F^{'}(x)+1=1+\varepsilon \ {\rm and} \
[x+F(x)]\varphi^{''}(x)+F^{''}(x)=-\frac{\varepsilon}{x}$$ from
which it follows that
$$a+\lambda =\frac{(1+\varepsilon)(1-\log\|g\|^2)}{\varepsilon},\ \
\frac{b}{a(b-1)}\leq
\frac{q\ell+\varepsilon+1}{\varepsilon+1}.$$Hence
the desired result follows directly from theorem 2.Q.E.D.\\

{\bf Remark 2}. Under the assumption that $\|g\|<1$ on $\Omega$, the
integrability
condition (21) is obviously weaker than (17).\\

We know by definition that $(0,-\frac{1}{2}
e^{-\varepsilon(x-1)},q)$ is another example of Skoda triples where
$\varepsilon$ is a positive constant and $q$ is the constant in
lemma 1. Thus theorem 2 applied to $(0,-\frac{1}{2}
e^{\varepsilon(x-1)},q)$ gives the following result.\\

{\bf Corollary 3} Let $\Omega\subseteq\mathbb{C}^n$ be a
pseudoconvex domain, $g_i\in\mathcal{O}(\Omega)(1\leq i\leq p)$,
$\psi\in{\rm PSH}(\Omega)$. We assume that $$\|g\|<1\ \text{holds on
}\Omega.$$ For every $f\in
{\bigwedge}^{\ell-1}\mathcal{O}(\Omega)^{\oplus p}$, if $\iota_gf=0$
and
$$\int_\Omega\|f\|^2
\|g\|^{-2(q\ell+1)} e^{-\psi}dV<\infty,\eqno{(23)}$$then there
exists some $u\in {\bigwedge}^{\ell}\mathcal{O}(\Omega)^{\oplus p}$
such that $\iota_gu=f$ and
$$\int_\Omega\|u\|^2 \|g\|^{2(-q\ell+\varepsilon)}e^{-\psi}dV
\leq C_\varepsilon\int_\Omega\|f\|^2\|g\|^{-2(q\ell+1)}
e^{-\psi}dV\eqno{(24)}$$where $p\in\mathbb{N}, 1\leq\ell\leq p,$
$\varepsilon$ and $C_\varepsilon$ are both positive constants($C_\varepsilon$ is
determined by $\varepsilon$) and $q$ is the constant in lemma 1.\\

Proof. By direct computations, we obtain
$$[x+F(x)]\varphi^{'}(x)+F^{'}(x)+1=\frac{\varepsilon}{2}e^{-\varepsilon(x-1)}+1, \ \
[x+F(x)]\varphi^{''}(x)+F^{''}(x)=-\frac{\varepsilon^2}{2}e^{-\varepsilon(x-1)}.$$
Hence we have for $\xi>1$
$$\frac{b}{a(b-1)}=\frac{1}{\xi-\frac{1}{2}
e^{-\varepsilon(\xi-1)}}+\frac{q\ell}{1+\frac{\varepsilon}{2}
e^{-\varepsilon(\xi-1)}} \leq 2+q\ell$$and
\begin{eqnarray*}a+\lambda &=& \xi-\frac{1}{2}e^{-\varepsilon(\xi-1)}
+2\varepsilon^{-2}(1+\frac{\varepsilon}{2}e^{-\varepsilon(\xi-1)})^2
e^{\varepsilon(\xi-1)}\\
&\leq&[\varepsilon^{-1}e^{\varepsilon
-1}+2(\frac{1}{\varepsilon}+\frac{1}{2})^2]e^{\varepsilon(\xi -1)}\\
&=& D_\varepsilon e^{\varepsilon(\xi-1)}\end{eqnarray*}where
$D_\varepsilon$ is a positive constant determined by $\varepsilon$.
Now corollary 3 follows from theorem 2 by choosing
$(\varphi,F,q)=(0,-\frac{1}{2}
e^{\varepsilon(x-1)},q)$ and the constant $C_\varepsilon$ in (24) could be taken to be $(2+q\ell)D_\varepsilon$. Q.E.D.\\

{\bf Remark 3.} (i) It is easy to see that when $\|g\|<1$ is valid
on $\Omega$ the integrability condition (23) in corollary 3 is
weaker than (21) but the estimate (22) for the solution in corollary
2 is stronger than (24). (ii) Comparing corollary 3 with theorem 1,
we see that if $\|g\|<1$ holds on $\Omega$ then the constant $\tau$
in theorem 1 could be chosen to be 1 (the coefficient
$\frac{\tau}{\tau-1}$ on the right hand of (18) should be replaced
by $C_\varepsilon$). (iii) It is interesting to compare corollary 3
with the main result of [T00] by setting $\ell=1,p=n$ and $
g_i=z_i(1\leq i\leq n).$ (iv) We may also choose the Skoda triple
more generally
 to be $(0,-\eta e^{-\varepsilon(x-1)},q)$ where $0<\eta<1$ is a constant,
but such a choice only results in a different constant
$C_\varepsilon $. (v) We can use the Skoda triple
$(\varepsilon_1\log x,-\varepsilon_2 e^{-\varepsilon_3(x-1)},q)$ to
combine the results in corollaries 2 and 3.
Here, $\varepsilon_1\geq 0,1>\varepsilon_2\geq 0,\varepsilon_3>0$ are constants satisfying $\varepsilon_1+\varepsilon_2>0.$\\

{\bf Final Comments.} As mentioned before, one can use the technique
of denominators to produce Skoda triples of the type $(0,F,q)$.
Hence we can deduce from our theorem 2 numerous examples of division
theorems. Actually, we can formulate a division theorem for the
Koszul complex in the same manner of theorem 2.7 in [V08]. To prove
this result, we just need to replace theorem 2.1 in [V08] by our
theorem 2 and then repeat its proof.\\

\section{ Appendix: Proof of Lemma 1}.

Let $V,W$ be Hermitian spaces with $\dim_\mathbb{C} V=p, \dim
_\mathbb{C}W=n,$ and $\{v_1,\cdots,v_p\},$ $\{w_1,\cdots,w_n\}$ be
orthonormal bases of $V,W$ respectively. We denote the dual bases by
$\{v^*_1,\cdots,v^*_p\}\subseteq V^*,\{w^*_1,\cdots,w^*_n\}\subseteq
W^*.$

Set\begin{eqnarray*}\mathcal {A}&=&\sum_{\overset{1\leq i\leq
p}{1\leq \alpha\leq n}}c_{ii_1\cdots i_{\ell-1}\alpha}w^*_\alpha\otimes v_i\in {\rm Hom}_\mathbb{C}(W,V),\\
\mathcal {B}_1 &=& \sum_{\overset{1\leq i\leq p}{1\leq \alpha\leq
n}}b_{i\alpha}v^*_i\otimes w_\alpha\in {\rm Hom}_\mathbb{C}(V,W),\\
X &=& \sum_{1\leq i\leq p}\overline{a_i}v_i\in V, \ \ \theta =
\sum_{1\leq i\leq p}a_iv^*_i\in V^*,\\ \mathcal {B} &=&
\iota_X(\theta\wedge\mathcal {B}_1)\in {\rm
Hom}_\mathbb{C}(V,W),\end{eqnarray*} then we know by definition the
following facts
$$\mathcal {A}\mathcal {B}=\iota_X(\theta\wedge\mathcal {A}\mathcal
{B}_1)\in {\rm End}_\mathbb{C}V$$ and $$ {\rm L.H.S. \ of }\
(8)=|{\rm Tr}\mathcal {A}\mathcal {B}|^2,\ \ \ \ {\rm R.H.S. \ of }\
(8)=q\|X\|^2\|\theta\wedge\mathcal {A}\mathcal {B}_1\|^2$$ where
both trace and norm are taken with respect to the Hermitian
structure on $V$. It remains therefore to show $$|{\rm Tr}\mathcal
{A}\mathcal {B}|^2\leq q\|X\|^2\|\theta\wedge\mathcal {A}\mathcal
{B}_1\|^2.$$Since the Cauchy-Schwarz inequality gives $$\|\mathcal
{A}\mathcal {B}\|^2\leq \|X\|^2\|\theta\wedge\mathcal {A}\mathcal
{B}_1\|^2$$ and $$|{\rm Tr}\mathcal {A}\mathcal {B}|^2\leq{\rm
Rank}_\mathbb{C}\mathcal {A}\mathcal {B}\|\mathcal {A}\mathcal
{B}\|^2,$$it suffices to estimate the upper bound of ${\rm
Rank}_\mathbb{C}\mathcal {A}\mathcal {B}.$

Since $v_{i_1\cdots i_\ell}$ are skew-symmetric on
$i_1,\cdots,i_\ell$, we get $${\rm Im} \mathcal {A}\mathcal
{B}\subseteq {\rm span}_{\mathbb{C}} \{
v_{i_1},\cdots,v_{i_{\ell-1}}\}^\bot.$$ On the other hand, we also
have
$$X\in {\rm Ker} \mathcal {A}\mathcal {B}.$$ We assume, without loss of
generality, $X\neq 0$ then we obtain the following estimate
\begin{equation*}{\rm Rank}_\mathbb{C} \mathcal {A}\mathcal
{B}\leq\begin{cases} \min\{p-1,n\}, & \ell=1;\\ \min\{p-\ell+1,n\},
& \ell\geq 2. \end{cases}\end{equation*} This is the desired rank estimate. Q.E.D.\\

Department of Mathematics

Fudan University

Shanghai 200433, China

E-mail address: qingchunji@fudan.edu.cn

\end{document}